\documentclass{amsart}          
\usepackage{amsmath,amsthm}     
\usepackage{amssymb}            
\usepackage{euscript}           
\usepackage{enumerate,calc}
\usepackage[matrix,arrow,curve,frame]{xy}    
\xymatrixcolsep{1.9pc}                          
\xymatrixrowsep{1.9pc}
\newdir{ >}{{}*!/-5pt/\dir{>}}                  
\numberwithin{equation}{section}

\newtheorem{thm}[equation]{Theorem}
\newtheorem{defn}[equation]{Definition}
\newtheorem{prop}[equation]{Proposition}

\newtheorem{cor}[equation]{Corollary}
\newtheorem{lemma}[equation]{Lemma}

\theoremstyle{definition}  
\newtheorem{example}[equation]{Example}

\newtheorem{note}[equation]{Note}

\newtheorem{remark}[equation]{Remark}
\newtheorem{reference}[equation]{Reference}

  {\end{list}}

\newcommand{\tens}              {\otimes}               

\newcommand{\cat}{\EuScript}    
\newcommand{\cC}{{\cat C}}

\newcommand{\cE}{{\cat E}}
\newcommand{\cF}{{\cat F}}

\newcommand{\cM}{{\cat M}}

\newcommand{\cR}{{\cat R}}

\newcommand{\Gr}{{\cat Grpd}}

\newcommand{\Set}{{\cat Set}}
\newcommand{\cCat}{{\cat Cat}}
\newcommand{\sSet}{s{\cat Set}}

\newcommand{\Pre}{Pre}

\newcommand{\field}[1]  {\mathbb #1} 

\newcommand{\Z}         {\field Z}

\DeclareMathOperator*{\colim}{colim}
\DeclareMathOperator*{\holim}{holim}
\DeclareMathOperator*{\hocolim}{hocolim}

\DeclareMathOperator{\Hom}{Hom}

\DeclareMathOperator{\ob}{ob}

\DeclareMathOperator{\sk}{sk}

\DeclareMathOperator{\eq}{eq}

\DeclareMathOperator{\Cat}{Cat}

\newcommand{\ra}{\rightarrow}                   
\newcommand{\lra}{\longrightarrow}              
\newcommand{\llra}[1]{\stackrel{#1}{\lra}}      

\newcommand{\we}{\llra{\sim}}                   

\newcommand{\blank}{-}                          

\newcommand{\m}{\mathcal M}         

\newcommand{\tuborg}{\left\{\begin{array}{ll}}
\newcommand{\sluttuborg}{\end{array}\right.}

\newcommand{\fpi}{{\pi_{oid}}}

\DeclareMathOperator{\Aut}{Aut}

\DeclareMathOperator{\Tot}{Tot}

\begin{document}

\email{sjh@math.huji.ac.il+}

\title{A Homotopy Theory for Stacks}

\author{Sharon Hollander}

\address{Department of Mathematics,
Hebrew University, Jerusalem, Israel}

\date{August, 1, 2005}

\begin{abstract}
We give a homotopy theoretic characterization of stacks on a site $\cC$
as the {\it homotopy sheaves} of groupoids on $\cC$.
We use this characterization to construct a model category
in which stacks are the fibrant objects.
We compare different definitions of stacks and show 
that they lead to Quillen equivalent model categories.
In addition, we show that 
these model structures are Quillen equivalent to the 
$S^2$-nullification 
of Jardine's model structure on sheaves of simplicial sets on $\cC$.

\end{abstract}

\maketitle

\section{Introduction}

The main purpose of this paper is to show that the classical definition 
of stacks \cite[Definition 4.1]{DM} ,\cite[II.$2.1$]{Gi} , can be stated in 
terms of homotopy theory.
From this point of view the definition appears natural, 
and places stacks into a larger homotopy theoretic context.
Constructions that are commonly performed on
stacks such as 2-category pullbacks, stackification, sheaves over a stack
and others, have easy homotopy-theoretic interpretations.

The basic idea is this: a stack $\m$ is
an assignment to each scheme $X$ of a groupoid $\m(X)$,
which is required to satisfy `descent conditions' \cite[II.$1.1$]{Gi} .
The descent conditions describe the circumstances under which
we require that local data glue together to yield global data.

Naively, one might propose a ``local-to-global'' requirement that 
the assignment ``isomorphism classes of -'' satisfy the sheaf condition.  
However, for very fundamental reasons, this almost never happens in examples.
Taking isomorphism classes is a localization process,
and such processes rarely preserve limits such as those which arise
in the statement of the sheaf condition.

Instead, one can ask that the assignment of groupoids satisfy
a sheaf condition with respect to the {\it best functorial approximation
to the limit} which respects isomorphism classes. This
is called the homotopy limit, denoted $\holim$, and is the total right
derived functor of inverse limit, (see \cite[Sections 9-10]{DS}). 

We prove that this {\it homotopy sheaf condition}
is exactly the content of the descent conditions and 
so, in this sense, stacks are the {\it homotopy sheaves}.

\begin{thm} \label{first-theorem}
Let $\cC$ be a Grothendieck topology.  
Let $\Gr/\cC$ denote the category of categories fibered in groupoids 
over $\cC$ (see Definition \ref{defFG}).
For each $\cE\in \Gr/\cC$, and $X\in \cC$, denote by $\cE(X)$ the groupoid of maps 
$\Hom(\cC/X,\cE)$.  
Then $\cE$ is a stack \cite[Definition $4.1$]{DM} if and only if
for every cover $\{U_i \ra X\}$ in $\cC$ we have an equivalence of categories
$$\xymatrix{ \cE(X) \ar[rr]^-{\sim} & & \holim \biggl(\prod \cE(U_i) 
\ar@2{->}[r] &
\prod \cE(U_{ij}) \ar@3{->}[r] & \prod \cE(U_{ijk}) \dots \biggr). }$$
Here $U_{i_0\dots i_n}$ denotes the iterated fiber product
$U_{i_0} \times_X \cdots \times_X U_{i_n}$, and we take the
homotopy limit of this cosimplicial diagram in the category of groupoids 
(see Section $2$).
\end{thm}

Classically stacks were defined either as categories fibered in groupoids or
lax presheaves \cite[Section 1]{Brn} satisfying descent conditions.
It is not hard to check that categories fibered in groupoids and
lax presheaves are equivalent categories and that these two definitions
agree \cite[Appendix B]{Holl}.
The following results make precise the sense in which it suffices to work with actual 
presheaves of groupoids instead of lax presheaves or categories fibered in groupoids.
Let $P(\cC,\Gr)$ denote the categories of presheaves of groupoids on $\cC$.

\begin{thm} \label{basic-mc}
There exists a model category structure on $\Gr/\cC$ in which weak 
equivalences are the fiberwise equivalences of groupoids.  Similarly, there exists model category 
structure on $P(\cC, \Gr)$ in which weak equivalences are the objectwise equivalences 
of groupoids.  For these two model structures the adjoint pair (see Section $3.3$)
 $$\xymatrix{ \Gr/\cC \ar@/^/[r]^{\Gamma}  & P(\cC,\Gr) \ar@/^/[l]^p }$$
is a Quillen equivalence.
\end{thm}

This model structure on $\Gr/\cC$ encodes many of the classical constructions for 
categories fibered in groupoids and for stacks. The equivalences of
stacks are just the weak equivalences. The $2$-category pullback is the homotopy pullback.
Writing these constructions in terms of standard homotopy theory makes clear their
functoriality properties and relations with one another.
 
This adjunction sends the full subcategory of stacks in $\Gr/\cC$ to the full subcategory
of $P(\cC,\Gr)$ of those presheaves which satisfy descent, which we call stacks. 

\begin{defn}
\label{defn-stacks-P}
Let $\cC$ be a Grothendieck topology. A presheaf of groupoids, $F$ on $\cC$
is a {\bf stack} if for every cover $\{U_i \ra X\}$ in $\cC$,
we have an equivalence of categories
$$\xymatrix{ F(X) \ar[rr]^-{\sim} & & \holim \biggl(\prod F(U_i) 
\ar@2{->}[r] &
\prod F(U_{ij}) \ar@3{->}[r] & \prod F(U_{ijk}) \dots \biggr). }$$
\end{defn}

The fact that $(\Gamma, p)$ is a Quillen equivalence tells us that it is equivalent 
to make any homotopy theoretic construction in one of these categories or the other,
(i.e. we can use $\Gamma$ to translate our problem in $\Gr/\cC$ into one about presheaves,
solve it and then apply $p$ to the result). 

Finally, we show that the model category structures of Theorem \ref{basic-mc}
can be localized and the result is a 
model structure even better suited to the study of stacks.
One can also consider the category of sheaves of groupoids on $\cC$, 
denoted $Sh(\cC,\Gr)$, instead of presheaves.

\begin{thm} \label{maintheorem}
There are simplicial model category structures $\Gr/\cC_L$, 
$P(\cC,\Gr)_L$, and $Sh(\cC,\Gr)_L$,
in which the stacks are the fibrant objects.\\
The adjoint pairs  
$$\xymatrix{ 
\Gr/\cC_L \ar@/^/[r]^{\Gamma}  & P(\cC,\Gr)_L \ar@/^/[l]^p \ar@/_/[r]_{sh}
& Sh(\cC,\Gr)_L, \ar@/_/[l]_{i}}$$
are Quillen equivalences (where the right adjoints point to presheaves).
All of these functors take stacks as defined in the domain category to 
stacks as defined in the range category and thus restrict to give
adjoint pairs between the stacks in each of these categories. 
\end{thm}

\begin{cor} \label{cor-stack}
In each of the above model categories the 
fibrant replacement functor gives a stackification functor.
\end{cor}




Of the different categories mentioned above, the simplest to analyze is
presheaves of groupoids which is closely related to 
presheaves of simplicial sets.  
Here a form of Dugger's local lifting conditions 
\cite[Section $3$]{DHI}, modified for groupoids, 
(see Definition \ref{local-equivalence-conditions})
allows us to characterize weak equivalences in a simple way (Theorem \ref{llc-thm}).
The comparison with the homotopy theory of simplicial sets is encapsulated 
in the following result.   

\begin{prop} \label{prop-jardine}
The above model structure on $P(\cC,\Gr)$ is Quillen equivalent to
Joyal's model category structure on $P(\cC,s\Set)$ localized
with respect to the maps $\partial \Delta^n \tens X \ra \Delta^n \tens X$,
for each $X \in \cC$ and  $n>2$.
\end{prop}

This theorem tells us that the homotopy theory of stacks is recovered from
Joyal's model category by eliminating all higher homotopies.
A direct corollary of this Quillen equivalence is the following result which 
nicely generalizes the usual criterion for a map to be an isomorphism between 
two sheaves of sets.

\begin{cor} \label{cor-ep}
If the topology on $\cC$ has enough points \cite[p.521]{MM}, the weak 
equivalences in $P(\cC,\Gr)$ are the stalkwise equivalences of groupoids.
\end{cor}

The characterization of stacks as the {\it homotopy sheaves} of groupoids
is the source and inspiration for all of the above results.
Furthermore this characterization generalizes naturally to a definition of
$n$-stack as follows: 

\begin{defn}  \label{definition-n-stack}
A presheaf of simplicial sets $F$ on $\cC$
is an {\bf $n$-stack} if for every $X$ in $\cC$,
$F(X)$ is a Kan complex, and
for every  hypercover $U_\bullet \ra X$ in $\cC$ \cite[Definition $4.2$]{DHI},
we have an equivalence of categories
$$\xymatrix{ F(X) \ar[rr]^-{\sim} & & \holim F(U_\bullet)} $$
where the homotopy limit is taken in the category of $(S^{n+1})^{-1}s\Set$, 
the $S^{n+1}$ nullification of simplicial sets.
\end{defn}

The reason for considering hypercovers as opposed to covers is discussed in
the introduction of \cite{DHI}.  There a model category is presented for 
presheaves of simplicial sets on $\cC$ in which these $n$-stacks are the fibrant 
objects, and the analogues of \ref{prop-jardine} and \ref{cor-ep} are proven.

For $n=1$, Definition \ref{definition-n-stack} is equivalent to \ref{defn-stacks-P},
by \cite[Theorem $1.1$]{DHI}, and Theorems \ref{maintheorem} and \ref{llc-thm}.
This is the foundation for the recent work of Toen and Vezzosi
on {\it homotopical algebraic geometry} (see \cite[Section $3$]{TV}),
following Simpson who first showed that model categories are useful
in understanding the theory of higher stacks, see \cite{HS}.

\subsection{Relation to other work}

In \cite{JT}, Joyal and Tierney introduce a model structure on
sheaves of groupoids on a site where the fibrant objects satisfy a 
strengthening of the stack condition, and are called strong stacks. 
It follows from Proposition \ref{comparison} that 
our model structures for stacks are Quillen equivalent to Joyal and Tierney's.

The main difference of our treatment is that we show that
the {\it descent conditions}, and hence
{\bf \emph{the classical definition of stack} } 
can be described in terms of a natural homotopy theoretic
generalization of the sheaf condition. It is this characterization which 
leads to a model category structure where the fibrant objects are 
\emph{precisely} the stacks.
In addition, our construction of the model structure draws a precise 
connection between the classical theory of stacks and the homotopy theory 
of simplicial presheaves.

\subsection{Contents}

In section $2$ we define the model structure on groupoids, and prove that
it is Quillen equivalent to a localization of simplicial sets with respect
to the map $S^2 \ra \ast$, called the $S^2$ nullification of $\sSet$.
We then present formulas for homotopy limits and colimits in $\Gr$
and prove that the descent category is a model for the homotopy inverse
limit of a cosimplicial diagram of groupoids.

In section $3$ we review the definition of categories fibered in groupoids 
over a fixed base category $\cC$, denoted $\Gr/\cC$.  
We construct an adjoint pair of functors between $\Gr/\cC$ and the category of 
presheaves of groupoids on $\cC$.  We show that these functors send the
subcategory of stacks in $\Gr/\cC$ to the subcategory of stacks in $P(\cC, \Gr)$.
Using this adjoint pair we prove Theorem \ref{first-theorem}.

In section $4$ we put model structures on $\Gr/\cC$ and $P(\cC, \Gr)$. 
In the basic model structure on each of these categories the weak equivalences 
are defined to be objectwise (or fiberwise). We note that the adjoint pair 
$(p,\Gamma)$ between these categories is a Quillen equivalence. 
We also observe that these model structures can
be localized with respect to the \emph{local} equivalences
$\{ \holim U_{\bullet} \ra X \}$.  In these local
model structures the fibrant objects are the stacks, 
and the adjoint pair $(p,\Gamma)$ is still a Quillen equivalence.  
This proves most of Theorem \ref{maintheorem}.

In section $5$ we define the local lifting conditions and prove 
Proposition \ref{prop-jardine}.
Finally, we show that there is a local model structure on sheaves of groupoids
on $\cC$ for which $(sh,i)$ is a Quillen equivalence, 
which completes the proof of Theorem \ref{maintheorem}.

Appendix $A$ contains proof that limits and colimits exist in the
category $\Gr/\cC$ of categories fibered in groupoids, which is needed 
to show that one can put a model structure on $\Gr/\cC$. 

Appendix $B$ contains a proposition about pushouts of categories needed
for the proof of the left properness of the model structure on $\Gr$ and 
for that on $\Gr/\cC$.

\subsection{Notation and Assumptions}

For a general introduction to Grothendieck topologies, see \cite{Ta},
\cite{MM}.
So as not to run into set theoretic problems, we assume that the 
Grothendieck topology $\cC$ is a small category.  

We write $\Pre(\cC)$ for the category of presheaves of sets on $\cC$.
For $\{U_i \ra X\}$ a cover in $\cC$, and $F$ a presheaf on $\cC$, let
\begin{itemize}
	\item the $n+1$-fold product $U\times_X U \times_X \dots \times_X U$ denote the coproduct 
		$\coprod U_{i_0} \times_X \cdots \times_X U_{i_n}$, 
		of the representable functors in $\Pre(\cC)$,
		where the coproduct is taken over all multi-indices $(i_0, \dots, i_n)$.
	\item $U_{\bullet}$ denote the simplicial diagram in $\Pre(\cC)$ with
         	$(U_{\bullet})_n$ equal to the $n+1$-fold product $U \times_X \cdots\times_X U,$
		with face and boundary maps defined by the
		various projection and diagonal maps. This is called
                the \emph{nerve} of the cover $\{U_i \ra X\}$.
	\item  $F(U \times_X U \times_X \dots \times_X U)$ denote the product 
	         $\prod F( U_{i_0} \times_X  \dots \times_X U_{i_n} )$, and
	         $F(U_{\bullet})$ the cosimplicial diagram $\Hom_{\Pre(\cC)}(U_\bullet,F)$.
\end{itemize}

For model categories and their localizations our references are
\cite{DS}, \cite{Ho}, \cite{Hi}.

\subsection{Acknowledgments}
The results in this paper are part of my doctoral thesis \cite{Holl}, MIT (2001), 
written under the supervision of Mike Hopkins.  
A great many thanks are owed to Dan Dugger, Gustavo Granja, and Mike Hopkins,
for many helpful conversations and ideas, without which this paper would
not exist.

\section{Homotopy Limits and Colimits of Groupoids}
We discuss here a model structure on the category of groupoids, denoted $\Gr$
and prove some results about homotopy limits and colimits which
will enable us in the next section to interpret the descent conditions in a
homotopy theoretic manner and prove Theorem \ref{first-theorem}.
The proofs of the results in this section are not hard.  For more details
the reader is referred to \cite{Holl}.

\subsection{Model Structure on $\Gr$}
The model structure we discuss on $\Gr$ is derived from that on $s\Set$
via the adjoint pair $(\fpi, N)$, where $\fpi$ denotes the fundamental
groupoid and $N$ the nerve construction.  Recall, $\fpi X$ is the groupoid
with objects $X_0$, and morphisms freely generated by $X_1$, subject to the 
relations $s_0 a = id_a$ for each $a \in X_0$, and $d_2 x \circ d_0 x = d_1 x$ 
for each $x \in X_2$.  It follows that all morphisms in $\fpi X$ are isomorphisms.  

We will sometimes abuse notation and denote the groupoid
$\fpi(\Delta^i)$ by $\Delta^i$.  This is the groupoid with $i+1$
objects and unique isomorphisms between them.  Similarly, we will
sometimes denote $\fpi(\partial \Delta^i)$ by $\partial \Delta^i$.  Let
$BG$ denote the groupoid with one object whose automorphism group is
the group $G$.

Note that the morphisms $\partial \Delta^i \ra \Delta^i, i=0,1,2$, are
$$\emptyset \ra \star, \quad \text{                  }\{ \star, \star \} \ra \Delta^1,
\quad \text{                   } \Delta^2 \times (B\Z \ra \star) .$$

\begin{thm} \label{MCG}
There is a left proper, simplicial, cofibrantly
generated model category structure on $\Gr$ in which:
\begin{itemize}
	\item weak equivalences are functors which induce an equivalence
		of categories,
	\item fibrations are the functors with the right lifting property
		with respect to the map $\Delta^0 \ra \Delta^1$,
	\item cofibrations are functors which are injections on objects.
\end{itemize}
The generating trivial cofibration is the morphism $\Delta^0 \ra \Delta^1$,
and the generating cofibrations are the morphisms
$\partial \Delta^i \ra \Delta^i, i=0,1,2$.
\end{thm}

\begin{note}
In this model category structure all objects are both fibrant and cofibrant,
so all weak equivalences are homotopy equivalences.
\end{note}
This model category structure appears in \cite{An, Bo}, 
a detailed description and proof can be found in \cite{St}, section $6$. 
The fact that it is simplicial and cofibrantly generated
is easy to check, and the left properness follows from the fact that
all objects are cofibrant, see \ref{B-1}.

\begin{cor}
With this model structure on $\Gr$, the adjoint pair 
$\fpi:\sSet \leftrightarrow \Gr:N$ is a Quillen pair.
\end{cor}

The following are important observations which we will use freely.

\begin{lemma} {\label{lemma1}}
Let $G \llra{f} H$ be a map of groupoids.
The following are equivalent:
\begin{itemize}
	\item $f$ is a weak equivalence in $\Gr$
	\item $Nf$ is weak equivalence in $\sSet$
\end{itemize}
Similarly, the following are equivalent:
\begin{itemize}
	\item $f$ is a (trivial) fibration in $\Gr$.
	\item $Nf$ is a (trivial) fibration in $\sSet$.
	\item $f$ has the right lifting property with respect to
		$\Delta^0 \to \Delta^1$ (with respect to
		$\partial \Delta^n \to \Delta^n$ for $n=0,1,2$).
\end{itemize}
\end{lemma}

Consider the model structure on $\sSet$ which is the localization 
\cite[Definition $3.3.1.1$]{Hi} 
of the usual model structure with respect to the map
$\partial \Delta^3 \ra \Delta^3.$
We will call this the {\it $S^2$-nullification} of $\sSet$ 
\cite{DF}, denoted $(S^2)^{-1}\sSet$.
Notice that the maps $\partial \Delta^n \ra \Delta^n,n>2,$ 
are all weak equivalences in $(S^2)^{-1}\sSet$,
and so $(S^2)^{-1}\sSet$ is also the localization of $s\Set$ with respect 
to this set of maps.

\begin{lemma}
In $(S^2)^{-1}\sSet$, weak equivalences are
the maps which induce an isomorphism on $\pi_0$ and $\pi_1$ at all base 
points.
\end{lemma}

\begin{thm}
The adjoint pair
$\fpi:(S^2)^{-1}\sSet \leftrightarrow \Gr:N $
is a Quillen equivalence.
\end{thm}

\subsection{Homotopy Limits and Colimits}

In \cite[$18.1.2, 18.1.8, 18.5.3$]{Hi}, explicit constructions are given of 
homotopy limit and colimit functors in arbitrary simplicial model categories.
Explicit formulas for homotopy limits and colimits in simplicial sets
go back to \cite[Section XI.$4$]{BK}.  Here we will give simplified formulas for 
homotopy limits and colimits in case the simplicial structure on the category 
derives from an {\it enrichment with tensor and cotensor} (see \cite{Db}) 
over $\Gr$.

\begin{defn}
Let $M$ be a category enriched with tensor and cotensor over $s\Set$.
We say that the simplicial structure on $M$ {\emph derives from an 
enrichment over $\Gr$},  
if $M$ is enriched with tensor and cotensor over $\Gr$, and if
for all $A,B \in M, K \in s\Set$, there are natural isomorphisms
$$s\Set(A,B) \cong N( \Gr(A,B)), \quad A^K \cong A^{\fpi K}, 
\quad A \otimes K \cong A \otimes \fpi K.$$
compatible with the natural isomorphism for each pair of simplicial sets 
$\fpi(K \times K') \cong \fpi K \times \fpi K'$.
\end{defn}

Our main concern will be the homotopy limit of a cosimplicial diagram, 
and dually the homotopy colimit of a simplicial diagram.
Our simplified formula for the former will allow us in the next section to
interpret the descent conditions for stacks in a homotopy-theoretic manner.

Let $\cC$ be a simplicial model category, and $I$ a small category. 
The homotopy limit of an $I$-diagram $X$ in $\cC$ with each $X(i)$ fibrant
is the equalizer of the two natural maps
$$ \prod_i X(i)^{N(I/i)} \Rightarrow
	\prod_{j \llra{\alpha} i} X(i)^{ N(I/j)},$$
where $I/i$ denotes the category of objects over $i$.
Similarly, the homotopy colimit of an objectwise cofibrant
$I$-diagram $X$ is the coequalizer of the two maps
$$\coprod_{i \llra{\alpha} j}  X(i) \tens  N(j/I) \Rightarrow
        \coprod_i X(i) \tens N(i/I), $$
where $j/I$ denotes the category of objects under $j$.

For $Y$ a fibrant object and ${\bf X} \in \cC^I$ objectwise cofibrant,
these functors satisfy the equation
\begin{equation} \label{adj-hocolim-holim}
s\Set(\hocolim {\bf X},Y) \cong \holim s\Set({\bf X},Y).
\end{equation}

\begin{thm}
Let $\cC$ be a simplicial model category whose simplicial structure
derives from an enrichment over $\Gr$, and let $X^{\bullet}$ be a cosimplicial
object in $\cC$, with each $X^i$ fibrant.  Then a model for
the homotopy inverse limit of $X^{\bullet}$
is given by the equalizer of the natural maps
$$\prod_{i=0}^2 (X^i)^{\Delta^i} \Rightarrow
	\prod_{[j]\to[i]}^{i\leq 2,j\leq 1} (X^i)^{\Delta^j}.$$
\end{thm}

\begin{proof}
First, notice that writing $\sk_2\Delta^\bullet$ for the $2$-skeleton of
$\Delta^{\bullet}$, the inclusion $\fpi \sk_2\Delta^\bullet \to \fpi
\Delta^{\bullet}$ is an isomorphism.
It follows that $\Tot X^\bullet$, the space of maps from
$\Delta^{\bullet}$ to $X^{\bullet}$,
is isomorphic to $\Tot_2 X^\bullet$, the space of maps from the
restriction $\Delta^\bullet|_{\Delta[2]}$ to $X^\bullet|_{\Delta[2]}$
where $\Delta[2]$ denotes the full subcategory of $\Delta$ with
objects $[0],[1]$ and $[2]$.
Since the map $\fpi \sk_1 \Delta^2 \lra \fpi\Delta^2$ is surjective,
$Tot_2 X^\bullet$ is given by the equalizer in the statement of the
theorem.

It now suffices to show that the homotopy limit of $X^\bullet$ is
naturally homotopy equivalent to $\Tot X^{\bullet}$.
Using the definition of the homotopy limit in a simplicial model
category given above, this is an easy consequence of the following 
Proposition.
\end{proof}

\begin{prop}
There is a homotopy equivalence of cosimplicial groupoids
$$F:\fpi N(\Delta/[\bullet]) \leftrightarrow \fpi \Delta^{\bullet}:G.$$
\end{prop}

\begin{proof}
Morphisms in $\fpi N(\Delta/[n])$ are generated by the commutative triangles
$$\xymatrix{ [k] \ar[rr] \ar[dr] & & [m] \ar[dl] \\ & [n] }$$
and their formal inverses.
Since the spaces $N(\Delta/[n])$ and $\Delta^n$ are contractible
in their fundamental groupoids there is a unique isomorphism between any
two objects. Therefore when defining functors between them it suffices
to set the values on objects. 

Let $\fpi N(\Delta/[n]) \llra{F_n} \fpi \Delta^n$ be the functor which sends
the object $[m] \ra [n]$ to the vertex $[0] \llra{e_m} [m] \ra [n]$, 
		where $e_k:[0] \to [k]$ sends $0$ to $k$.
One can check easily that $F$ is natural in $n$,
and so defines a morphism
$$\fpi N(\Delta/[\bullet]) \llra{F} \fpi \Delta^{\bullet} \in c\Gr.$$
Let $G_n$ be the functor which is defined
on objects by including $[0] \ra [n]$ in $\Delta/[n]$.
Again it is easy to check that $G_n$ is natural in $[n]$,
and so defines a morphism
$\fpi \Delta^{\bullet} \llra{G} \fpi N(\Delta/[\bullet])$ in $c\Gr$.

The composition $F \circ G$ is the identity.
There are unique natural transformations $ G_n \circ F_n \llra{\alpha_n} id$
which must commute with the simplicial operations.
\end{proof}

The groupoid $\Tot_2(X^{\bullet})$ will also be called the
{\it descent category} of $X^{\bullet}$. From now
on, when we refer to the homotopy limit of a cosimplicial
groupoid $X^\bullet$ we will mean the simpler model $\Tot_2(X^{\bullet})$.
The following corollary gives an explicit description of this
groupoid.

\begin{cor} \label{desc-formula}
The homotopy inverse limit of a cosimplicial groupoid $X^{\bullet}$
is the groupoid whose
\begin{itemize}
	\item objects are pairs $(a,d^1(a)\llra{\alpha} d^0(a))$, with
		$a \in obX^0,\alpha \in morX^1$,
		such that $s^0(\alpha)=id_a$, and
		$d^0(\alpha)\circ d^2(\alpha)=d^1(\alpha)$,
	\item morphisms $(a,\alpha) \ra (a',\alpha')$ are maps
		$a \llra{\beta} a'$, such that the following diagram commutes
	$$\xymatrix{
		d^1(a) \ar[r]^{d^1(\beta)} \ar[d]^{\alpha} & d^1(a')
		\ar[d]^{\alpha'} \\ d^0(a) \ar[r]^{d^0(\beta)} & d^0(a')}$$
\end{itemize}
\end{cor}

Dually we have the following theorem giving a formula for
homotopy colimits of simplicial diagrams.

\begin{thm}
Let $\cC$ be a simplicial model category whose simplicial structure
derives from a groupoid action and let $X_{\bullet} \in s\cC$, be such that
each $X_i$ is cofibrant. Then the homotopy colimit of $X_{\bullet}$
is naturally homotopy equivalent to the coequalizer of the maps
$$\coprod_{[n]\lra [m]}^{n \leq1,m \leq 2} X_m \tens \Delta^n
        \Rightarrow \coprod_{n=0}^2 X_n \tens \Delta^n .$$
\end{thm}

\section{Stacks}
There are different categories in which
the descent condition can be formulated and so in which stacks can be defined.
In this section we will discuss stacks in the context of 
{\it categories fibered in groupoids} over $\cC$, denoted $\Gr/\cC$
\cite{DM,Gi}.

After discussing some important properties of  $\Gr/\cC$,
we will construct an adjoint pair
$ p : P(\cC,\Gr) \leftrightarrow \Gr/\cC: \Gamma , $
satisfying the following properties:
\begin{itemize}
	\item For $F$ in $P(\cC,\Gr)$, the map $F(X) \ra \Gamma pF(X)$
		is an equivalence of groupoids, for all $X \in \cC$.
	\item For $\cE \in \Gr/\cC$, the map $ p \Gamma \cE \ra \cE$
		is an equivalence of categories over $\cC$.
\end{itemize}

In the following subsection we will discuss the classical definition of stacks 
in $\Gr/\cC$ used in algebraic geometry \cite{DM}. 
We show that it can be reformulated in terms of homotopy limits of groupoids,
expressing stacks as those objects satisfying the
{\emph homotopy sheaf condition}.

\subsection{Categories Fibered in Groupoids over $\cC$}
\begin{defn} \cite{DM} \label{defFG}
The category $\Gr/\cC$ is the full subcategory of $\cCat/\cC$
whose objects are functors $\cE \llra{F} \cC$ satisfying the following
properties:
\begin{enumerate}
	\item Given $Y\llra{f} X \in \cC$, and $X'\in \cE$ such that
              $F(X')=X$, there exists $Y' \llra{f'} X' \in \cE$ such that
	      $F(f')=f.$
	\item Given a diagram in $\cE$, over the commutative diagram in $\cC$,
	$$\xymatrix{& Y' \ar[d]^{f'} & \ar@2{->}[r]^F & &
		&  Y \ar[dl]_h \ar[d]^f \\
		Z' \ar[r]^{g'}  & X'  & \ar@2{->}[r]^F & & Z \ar[r]^g & X, }$$
		with $F(f')=f, F(g')=g$, there exists a unique $h'$
		such that $g' \circ h' = f'$ and $F(h')=h$.
\end{enumerate}
\end{defn}

This definition may seem involved but it becomes very simple when we
look at the functors $F_{X'}$ induced by $F$ on the over categories
$$\cE/X' \llra{F_{X'}} \cC/X,$$
where $X' \in \cE$, and $F(X')=X$.
The conditions for $\cE \llra{F} \cC$ to be a category
fibered in groupoids over $\cC$ are equivalent to the
simple requirement that each functor $F_{X'}$ induce a surjective 
equivalence of categories.

Let $\cE_X$ denote the fiber category over $X$.  This
has objects those of $\cE$ lying over $X$
and morphisms those of $\cE$ lying over $id_X$.
It is easy to see that if $\cE \ra \cC \in \Gr/\cC$,
the fiber categories $\cE_X$ are groupoids.

\begin{example}
The simplest examples of categories fibered in groupoids over $\cC$
are the projection functors $\cC/X \ra \cC$ for each $X \in \cC$.
If $Y \llra{f} X$ is an object of $\cC/X$, then $(\cC/X)/f \cong \cC/Y$,
and so conditions $1.$ and $2.$ above are trivially satisfied.
Notice that $(\cC/X)_Y$ is the discrete groupoid whose set of objects
is $\Hom_{\cC}(Y,X)$.

Another class of simple examples are $\cC \times G \llra{pr} \cC$,
for $G \in \Gr$.
Here the fibers over each $X \in \cC$ are canonically isomorphic to $G$.
\end{example}

The proofs of the following results are straightforward.
\begin{lemma}
$\Gr/\cC$ is enriched with tensor and cotensor over $\Gr$.  The
objects of $\Gr(\cE,\cE')$ are the functors $\cE \ra \cE'$ over $\cC$,
and the morphisms are the natural isomorphisms between such functors
covering the identity natural automorphism of $id_{\cC}$.
Moreover, the tensor is given by the formula
$$\cE \tens G := \cE \times_{\cC} (\cC \times G),$$
and the cotensor $\cE^G$ is the category of functors from
$(G \ra *)$ to $(\cE \ra \cC)$.
\end{lemma}

\begin{prop}
Let $\cE \llra{F} \cC \in \Gr/\cC$, and $X \in \cC$, then
\begin{enumerate}
	\item for any $X' \in \cE$ with $F(X')=X$ there is a section
		$$\xymatrix{ & \cE \ar[d]^F \\ \cC/X \ar@{-->}[ur]^{G} \ar[r] &
			 \cC}$$ such that $G(id_X)=X'$.
	\item If $G,G':\cC/X \ra \cE$ are two sections and
		$G(id_X) \llra{f} G'(id_X)$ is a morphism $\cE_X$,
		then there is a unique natural isomorphism $G \llra{\phi} G'$
		over $id_{\cC}$, with $\phi(id_X)=f$.
\end{enumerate} \label{rig}
It follows that, for each $X\in \cC$, the natural map
$$\Gr(\cC/X, \cE) \ra \cE_X$$
given by evaluation at $id_X$ is a surjective equivalence of groupoids.
There is a left inverse which is unique up to unique natural isomorphism.
\end{prop}

This says that given $\cE \ra \cC$ there is a functorial
``rigidification'' of the fibers.
Later we will use this method of rigidification
to construct a functor from $\Gr/\cC$ to $P(\cC,\Gr)$.

The following observation, which can be proven in a similar fashion,
will be used in the next subsection.

\begin{prop}
Let $\cE \ra \cC$ be a category fibered in groupoids,
and $Y \llra{f} X$ morphism in $\cC$.
There are ``pullback'' functors $\cE_X \llra{f^*} \cE_Y$
which are unique up to a unique natural isomorphism.
\end{prop}

\subsection{Stacks} 

Let $\cE \ra \cC$ be a category fibered in groupoids,
and assume that for each $X \llra{f} Y$ we have chosen pullback functors
$\cE_Y \llra{f^*} \cE_X$.  Given a morphism $U_i \ra U \in \cC$,
we will sometimes abuse notation and denote the pullback of an element
$a \in \cE_U$ to $\cE_{U_i}$ by $a|_{U_i}$.
In defining some of the maps below, we will also make implicit use
of the natural isomorphisms $(a |_{U_i}) |_{U_{ij}} \cong a|_{U_{ij}}$.

\begin{defn} \cite{Gi} \cite{DM} \label{defn-stacks-DM}
A {\it stack} in $\Gr/\cC$ is an object $\cE \ra \cC$
which satisfies the following properties for any cover $\{ U_i \ra X \}:$
\begin{enumerate}
	\item given $a,b \in \cE_X$, the following is an equalizer sequence
	     $$\Hom_{\cE_X}(a,b) \ra \prod \Hom_{\cE_{U_i}}(a|_{U_i},b|_{U_i})
	     \Rightarrow \prod \Hom_{\cE_{ U_{ij}}}(a|_{U_{ij}},b|_{U_{ij}}),$$
	\item given $a_i \in \cE_{ U_i }$ and isomorphisms
		$$ a_i |_{ U_{ij} } \llra{ \alpha_{ij} } a_j |_{ U_{ij},}$$
		satisfying the \emph{cocycle condition}
		$$\alpha_{jk} |_{ U_{ijk} } \circ \alpha_{ij} |_{ U_{ijk} } =
		\alpha_{ik} |_{ U_{ijk}, }$$
		then there exist $a \in \cE_X$, and isomorphisms
		$a|_{U_i} \llra{\beta_i} a_i$, such that the following square
		commutes
		\begin{equation} \label{square}
                        \xymatrix{ a|_{U_{ij}} \ar[r]^{\beta_i |_{ U_{ij} }}
			\ar[d]_= & a_i|_{U_{ij}} \ar[d]^{\alpha_{ij}} \\
			a|_{U_{ij}} \ar[r]^{\beta_j |_{U_{ij}}} &
			a_j |_{ U_{ij} } .}
                \end{equation}
\end{enumerate}
In this case, we say that $\cE \ra \cC$
\emph{satisfies descent}.
\end{defn}

\begin{note}\label{self-intersect}
The cocyle condition applied to indicies $(i,i,j)$ implies that
$\alpha_{ii} |_{ U_{iij} } = id$ which implies that $\alpha_{ii}$ is itself
the identity by $(1)$ of \ref{defn-stacks-DM}.
\end{note}

This definition seems very complicated, but it can be considerably
simplified if we recall the description of the homotopy inverse limit of a
cosimplicial groupoid given in Corollary \ref{desc-formula}.


\begin{thm}[Theorem $1.1$]

A category fibered in groupoids $\cE \ra \cC$ is a
stack if and only if for all covers $\{ U_i \ra X\}$
\begin{equation} \label{map}
\Gr(\cC/X,\cE) \ra \holim  \Gr(\cC/U_{\bullet},\cE)
\end{equation}
is an equivalence.
\end{thm}

\begin{proof}
We begin by showing that condition $(1)$ in Definition \ref{defn-stacks-DM}
is equivalent to the requirement that for objects
$F_a, F_b \in \Gr(\cC/X,\cE)$, the set of
morphisms $F_a \ra F_b$ is in bijective correspondence
with the set of morphisms between their images
in $\holim \Gr(\cC/U_{\bullet},\cE)$.

Consider objects $F_a, F_b \in \Gr(\cC/X,\cE)$, and let $a=F_a(id_X)$ and
$b=F_b(id_X)$ in $\cE_X$.  Evaluation at $id_{(-)}$ induces bijections
$$\xymatrix{\Hom(F_a, F_b) \ar[r] \ar[d]^{\cong} &
\prod \Hom(F_a|_{U_i},F_b|_{U_i}) \ar@2{->}[r] \ar[d]^{\cong} &
\prod \Hom(F_a|_{U_{ij}}, F_b|_{U_{ij}}) \ar[d]^{\cong} \\
\Hom_{\cE_X}(a,b) \ar[r] & \prod \Hom_{\cE_{U_i}}(a|_{U_i},b|_{U_i})
\ar@2{->}[r] & \prod \Hom_{\cE_{U_{ij}}}(a|_{U_{ij}},b|_{U_{ij}}). }$$
It follows that the top line is an equalizer if and only if the bottom one 
is.
By corollary \ref{desc-formula},
the top line is an equalizer if and only if
$\Hom(F_a, F_b)$ is in bijective correspondence with the set of
maps from the image of $F_a$ to the image of $F_b$ in
$\holim \Gr(\cC/U_{\bullet},\cE)$.
The requirement that the bottom line be an equalizer is
condition $(1)$ in Definition \ref{defn-stacks-DM}.

To finish the proof we have to show that condition $(2)$
is equivalent to the requirement that
every object in $\holim \Gr(\cC/U_{\bullet},\cE)$ be
isomorphic to one in the image of $\Gr(\cC/X,\cE)$.
This follows from the description of morphisms in Corollary 
\ref{desc-formula}
once we show that specifying an object in
$\holim\Gr(\cC/U_{\bullet},\cE)$
is equivalent to specifying descent datum as in condition $(2)$
of Definition \ref{defn-stacks-DM}.

By corollary \ref{desc-formula}, an object of
$\holim \Gr(\cC/U_{\bullet},\cE)$, consists of
an object $F_c \in \prod \Gr(\cC/U_i,\cE)$, and an isomorphism
$d^1 F_c \llra{\alpha} d^0 F_c$, satisfying
$d^0 (\alpha) \circ d^2 (\alpha) = d^1(\alpha)$ and $s^0(\alpha)=id_{F_c}$.
For any $U \llra{f} V$, and $F_a \in \Gr(\cC/V, \cE)$ with $F_a(id_V)=a$,
the evaluation $F_a|_U (id_U)$ is a choice of pullback of $a$ along $f$,
and so $F_a|_U (id_U)$ is canonically isomorphic to the pullback $f^*a$,
which we chose in advance.
Evaluating at $id_{U_i}$ determines
$c \in \prod \cE_{U_i}$, and isomorphisms
$\alpha_{ij}=\alpha(id_{U_{ij}})$ satisfying the cocycle condition.
Composing with the canonical isomorphisms
$c|_{U_{ij}} \cong F_c|_{U_{ij}} (id_{U_{ij}})$, we obtain
isomorphisms $c|_{U_i} \llra{\bar{\alpha}_{ij}} c|_{U_j}$, satisfying the
cocycle condition.

Conversely, given $c \in \prod \cE_{U_i}$ and $\alpha_{ij}$,
as in condition $(2)$ satisfying $\Delta^*(\alpha_{ii})=id_{U_i}$ (see
Note \ref{self-intersect}), we can lift them to an object
$F_c \in \prod \Gr(\cC/U_i,\cE)$, and an isomorphism
$d^1 F_c \llra{\alpha} d^0 F_c$.  Since these lifts are essentially unique
they must also satisfy the cocycle condition and $s^0(\alpha)=
id_{F_c}$ and hence determine an object of $\holim 
\Gr(\cC/U_{\bullet},\cE)$.
\end{proof}

\subsection{Adjoint Pair Between $\Gr/\cC$ and $P(\cC,\Gr)$}

Let $\cE \ra \cC$ be a category fibered in groupoids. By Corollary 
\ref{rig},
the assignment to each $X \in \cC$ of the sections
$\Gr(\cC/X,\cE)$ is a functor such that $\Gr(\cC/X,\cE) \we \cE_X$.

\begin{defn}
\label{rigfunctor}
Let $\Gamma:\Gr/\cC \ra P(\cC,\Gr)$ be the functor which sends
$\cE \ra \cC$ to the presheaf $\Gamma \cE(X) := \Gr_{\Gr/\cC}(\cC/X,\cE)$.

Let  $p:P(\cC,\Gr) \ra \Gr/\cC$ be the functor defined by setting $pF$
to be the category whose
\begin{itemize}
	\item objects are pairs $(X,a)$ with $a \in F(X)$,
	\item morphisms $(X,a) \ra (Y,b)$ are pairs $(f,\alpha)$ where
		$X \llra{f} Y \in \cC$ and $a \llra{\alpha} F(f)(b)$
		is an isomorphism in $F(X)$.
\end{itemize}
The composition of two morphisms
$(X,a) \llra{(f,\alpha)} (Y,b) \llra{(g,\beta)} (Z,c)$
is the pair $(g\circ f, F(f)(\beta) \circ \alpha)$.
\end{defn}

It is easy to check that both $p$ and $\Gamma$ preserve the 
enrichment over $\Gr$ with tensor and cotensor.

\begin{thm} \label{APpgamma}
The functors
$p:P(\cC,\Gr)\leftrightarrow \Gr/\cC:\Gamma$
form an adjoint pair. The unit of the
adjunction is an objectwise equivalence, and the counit is a fiberwise
equivalence of groupoids.
\end{thm}

This adjoint pair restricts to an adjunction between the subcategory
of stacks and $\Gr/\cC$ and the subcategory of presheaves $P \in P(\cC,\Gr)$
which satisfy the following condition:

\begin{itemize} 

\item 
For any cover $\{ U_i \to X \}$ in $\cC$ the induced
map $P(X) \ra \holim P(U_\bullet)$ is an equivalence of groupoids. 

\end{itemize}
With this motivation we will also call a presheaf of groupoids satisfying the 
above condition a stack, (Definition \ref{defn-stacks-P}).

\section{Model Structures}
In this section we put model structures on $P(\cC,\Gr)$
and $\Gr/\cC$.
We first construct basic model structures,
then we localize them so that the local  
weak equivalences detect the topology on $\cC$.
We then observe that in these {\it local} model structures, the
fibrant objects are the stacks, and the weak equivalences are the maps 
which locally are an equivalence of groupoids.
Isomorphisms between sheaves can be detected locally
and this property characterizes the subcategory of sheaves, 
It follows from our analysis that analogously global equivalences between 
stacks can be detected locally, and this property is a characterization of stacks.

We will use the notation of \cite[Definition $7.1.3, 9.1.5$]{Hi} for the
model category axioms.

\subsection{The Basic Model Category Structures}

Henceforth we will abuse notation and denote by $X$ the representable
functor $\Hom_\cC(\blank,X)$ considered as a discrete groupoid.

\begin{prop} \label{MC-P-Sh}
There is a left proper, cofibrantly generated, model category structures on
$P(\cC,\Gr)$, where
\begin{itemize}
	\item $f$ is a weak equivalence or a fibration if
		 $\Gr(X,f)$ is one for all $X \in \cC$,
	\item cofibrations are the maps with the left lifting property with
		respect to trivial fibrations.
\end{itemize}
The maps of the form
$X \ra X \tens \Delta^1$, for $X \in \cC$, form a set of generating
trivial cofibrations.
The maps of the form
$X \tens \partial \Delta^i \ra X \tens \Delta^i$
for $X \in \cC$ and $i=0,1,2$ form a set of generating cofibrations.
\end{prop}

\begin{proof}
The proof is an easy exercise. 

\end{proof}

Now we construct a model category on $\Gr/\cC$
relative using the set of generators $\cC /X \ra \cC$.

\begin{thm} \label{MC-FG}
There is a left proper, cofibrantly generated, simplicial model category
structure on $\Gr/\cC$ in which
\begin{itemize}
	\item $f$ is a weak equivalence or a fibration if
		$\Gr_{\Gr/\cC}(\cC/X,f)$ is one for all $X \in \cC$,
	\item cofibrations are the maps with the left lifting property with
		respect to trivial fibrations.
\end{itemize}
The maps of the form
$\cC/X  \ra (\cC/X \tens \Delta^1)$, for $X \in \cC$, form a set
of generating trivial cofibrations.
The maps of the form
$(\cC/X  \tens \partial \Delta^i) \ra (\cC/X \tens \Delta^i)$,
for $X \in \cC$ and $i=0,1,2$ form a set of generating cofibrations.
\end{thm}
\begin{proof}
For M1, see Appendix A. M2-M4(1) are obvious.
In order to apply the small object argument to prove M5, we need to check
that the objects $\cC/X \tens G \ra \cC$ with 
$G=(\partial)\Delta^i,i=0,1,2$,
are small with respect to the colimits which appear in the small object
argument. First notice that sequential colimits in
$\Gr/\cC$ agree with sequential colimits in $\Cat/\cC$.
For convenience, in the construction of the factorization for M5(1)
we will take pushouts along both the generating cofibrations and
the generating trivial cofibrations.

Let $\cE_i \to \cE_{i+1}$ be constructed as usual, using the small object
argument, and let consider a map $F:\cC/X \lra \colim \cE_i$.
$F(id_X)$ lifts to some element $X'$ in some $\cE_i$,
and we can extend this to a map $F_i':\cC/X \lra \cE_i$.
Let $F'$ be the composition $\cC/X \ra \cE_i \ra \colim \cE_i$.
Then $F'(id_X)=F(id_X)$, and so there is a unique natural
isomorphism $\phi:F\lra F'$ making the following diagram commute
$$\xymatrix{ & \cC/X  \ar[r]^{F_i'} \ar[d] \ar[dr]^{F'} & \cE_i \ar[d] \\
\cC/X \ar[r]  \ar@/_3ex/[rr]_{F} & \cC/X \tens \Delta^1 \ar@{-->}[r]^{\phi}
& \colim \cE_i.}$$
The map $\cC/X \ra \cC/X \tens \Delta^1$ is one of the
generating trivial cofibrations, so by construction we obtain a lift
$$\xymatrix{ & \cC/X \ar[r]^{F_i'} \ar[dd] & \cE_i \ar[d] \\ &  &
\cE_{i+1} \ar[d] \\  \cC/X \ar[r] \ar@{-->}[urr]  \ar@/_3ex/[rr]_{F} &
\cC/X \tens \Delta^1 \ar[r]^{\phi} \ar@{-->}[ur] & \colim \cE_i.}$$
Thus $\cC/X$ is small with respect $\colim \cE_i$. Since natural
transformations between sections are determined uniquely by their
evaluation on $id_X$, a similar argument shows that $\cC/X \tens
(\partial)\Delta^i$ is small with respect to $\colim \cE_i$.
This completes the proof of M5(1).

For M5(2) use the small object argument for the generating trivial cofibrations.
We need to show that the first map in the factorization is a weak equivalence.
Note that if $\cE \ra \cE'$ has the left lifting property with
respect to all fibrations, then in particular it has the left
lifting property with respect to $\cE \ra \cC$ and
$(\cE')^{\Delta^1} \ra (\cE')^{\partial \Delta^1}$, and therefore
it is an equivalence of categories over $\cC$. An equivalence of
categories over $\cC$ is clearly a weak equivalence. It follows that
the cofibration constructed using the small object argument for
M5(2) is also a weak equivalence.

M4(2) now follows by the same argument given in the
proof of Theorem \ref{MC-P-Sh}. M7 follows immediately from the
definition of (trivial) fibration in $\Gr/\cC$ and the adjunction
formulas given by the simplicial structure.

To show left properness, it suffices to show that the pushout of a
trivial fibration along a cofibration is a weak equivalence.
We begin by noting that trivial fibrations are surjective
equivalences of categories.
Let $F:\cE' \lra \cE''$ be a trivial fibration and let $X', Y' \in \cE',
X''=F(X'), Y''=F(Y')$.  Clearly $F$ is surjective on objects and morphisms.
We will show that the map
$$\Hom_{\cE'}(X', Y') \ra  \Hom_{\cE''}(X'', Y'')$$
is a bijection.  If $F(f')=F(g')$ then $f'$ and $g'$ have the same image in
$\cC$ and so there is a unique isomorphism $h'$
filling in the following triangle in $\cE'$:
$$\xymatrix{ X' \ar[dr]^{f'} \ar@{-->}[d]^{h'} \\ X' \ar[r]^{g'} & Y'.}$$
By the uniqueness of the lifting $h'$, $F(h')=id_{X''} \in \cE''$.
Since $F$ is a trivial fibration it follows that $h'=id_{X'}$.

Now note that cofibrations in $\Gr/\cC$ are inclusions on objects
as this is the case for the generating cofibrations.
Proposition \ref{B-1} implies that the pushout in $\Cat/\cC$ of
a surjective equivalence of categories along an inclusion on objects is
still an equivalence of categories over $\cC$. This simultaneously implies 
that
the pushout in $\Cat/\cC$ coincides in this case with the pushout in 
$\Gr/\cC$
(see the proof of Theorem \ref{complete}) and completes the proof.
\end{proof}
\begin{cor}
\label{quiequivcatpre}
The adjoint pair $p :P(\cC,\Gr) \leftrightarrow \Gr/\cC: \Gamma$
is a Quillen equivalence.
\end{cor}

\subsection{Local Model Category Structures}

For convenience, we will now also denote by $X$ the category fibered 
in groupoids $\cC/X \to \cC$. In the $P(\cC,\Gr)$
or $\Gr/\cC$, let $S$ denote the set of maps
\[ S=\{ \hocolim U_\bullet \to X : \{U_i \to X\}\text{ is a cover in } \cC\} 
\]
where $U_{\bullet}$ denotes the nerve of the covering
$\{U_i \to X\}$.

\begin{prop} \label{loc-mod-str}
Let $\cM$ be $P(\cC,\Gr)$ or
$\Gr/\cC$. There is a model category structure on $\cM$ which
is the localization of the model structure of Theorems \ref{MC-P-Sh} or
\ref{MC-FG} with respect to the set of maps $S$.
\end{prop}
We call these weak equivalences {\it local weak equivalences}.
\begin{proof}
Since homotopy colimits of cofibrant objects are cofibrant, the domains
and ranges of the morphisms in the localizing set are cofibrant.
By Theorems \ref{MC-P-Sh} and \ref{MC-FG}, the model category structures
on $P(\cC,\Gr)$ and $\Gr/\cC$ satisfy the hypothesis of
\cite[Theorem 4.1.1]{Hi}, so the proposition follows.
\end{proof}

Let $\cM$ be $P(\cC,\Gr)$ or
$\Gr/\cC$. We will write $\cM_L$ for the category $\cM$ with the model
structure given by the previous proposition.

\begin{cor}
The adjoint pair $p :P(\cC,\Gr)_L \leftrightarrow (\Gr/\cC)_L: \Gamma$
is a Quillen equivalence.
\end{cor}

Since in the old model structure on $\cM$ every object is fibrant,
and $X \in \cC$ is cofibrant,
an object $F \in \cM_L$ is fibrant if and only if
$$\Gr(X,F) \ra \Gr(\hocolim U_{\bullet} ,F) = \holim \Gr(U_{\bullet},F)$$
is a weak equivalence for all covers. This
happens if and only if  $F$ is a stack.
It follows that {\bf a fibrant replacement functor for $\cM_L$
is a stackification functor}.
One of the properties of localizations of model categories is that local 
equivalences between fibrant objects are just the old equivalences.  It follows
that {\bf a local equivalence between stacks is just an objectwise weak 
equivalence}.  

\begin{remark} Since stacks are the fibrant objects, and  representables
are cofibrant, it follows that when $\m$ is a stack,
$h\Hom(X,\m)$ is equivalent to the groupoid $\m(X)$.
In particular, $[X,\m]$ is the set of isomorphism classes of $\m(X)$.
\end{remark}

\begin{remark}
It is not hard to check that a small presentation (in the
sense of \cite[Definition 6.1]{Dg}) of $P(\cC,\Gr)_L$
is given by the Yoneda embedding of $\cC$ in $P(\cC,\Gr)$ and
the set of maps \begin{itemize}
\item $X \tens \partial \Delta^n \ra X \tens \Delta^n,$ 
  for all $X \in \cC, n>2$,
\item $\hocolim U_\bullet \ra X$ for all covers 
   $\{U_i \ra X\}$  in $\cC.$ \end{itemize}
This means that the local model category structure is the ``quotient''
of the universal model category generated by $\cC$ by the
relations given by the maps above.
\end{remark}

\section{Characterization of Local Equivalences}

In this section we prove that a morphism $f$ is a local weak equivalence in
the model structure of proposition \ref{loc-mod-str}
if and only if it satisfies one of the following equivalent properties:
\begin{itemize}
\item $f$ is an {\it isomorphism on sheaves of homotopy groups},
\item $f$ satisfies the {\it local lifting conditions}, 
		(Definition \ref{local-equivalence-conditions}),
\item $f$ is a {\it stalkwise weak equivalence} (when $\cC$ has
enough points).
\end{itemize}
It follows that a map between stacks satisfying one of the above properties
is actually an objectwise equivalence.

We use this to prove that our local model structure $P(\cC, \Gr)_L$
is Quillen equivalent to the $S^2$-nullification of Joyal's model
structure on presheaves of simplicial sets \cite{Ja}.

Using the characterization of local weak equivalences we also prove that there 
is a local model category structure on $Sh(\cC, \Gr)$ such that the adjoint 
pair $sh:P(\cC, \Gr)_L  \leftrightarrow Sh(\cC, \Gr):i$ 
is a Quillen equivalence.

\subsection{Joyal's Model Structure}
For a simplicial set $X$, and basepoint $a \in X_0$, $\pi_n(X,a)$ 
denotes the $n$-th homotopy group of the fibrant replacement of $X$ 
with basepoint the image of $a$.

\begin{defn} \cite{Ja}
Let $F$ be a presheaf of simplicial sets or groupoids. Then
\begin{itemize}
	\item $\pi_0 F$ is the presheaf of sets defined by
		$(\pi_0 F)(X):= \pi_0 (F(X))$. 
	\item For $F\in P(\cC,\sSet)$ and  $a \in F(X)_0$,
              $\pi_n(F,a)$ is the presheaf of groups
		on $\cC/X$ defined by
	$$\pi_n(F,a)(Y \llra{f} X) = \pi_n (F(Y),f^*a).$$
	For $F \in P(\cC,\Gr)$ and $a \in \ob F(X)$,
      $\pi_n(F,a):=\pi_n(NF,a)$.
\end{itemize}
We say that a map $F \llra{\phi} G$ of presheaves of simplicial sets or
groupoids is an \emph{isomorphism on sheaves of homotopy groups} if
the induced maps $sh \pi_0(\phi)$ and $sh \pi_n(\phi,a)$ are isomorphisms
for all $a \in F(X)$, and all $X \in \cC$.
\end{defn}
Note that if $F$ is a presheaf of groupoids then $\pi_i(F,a) = 0$
for $i>1$, and $\pi_1(F,a)$ is the presheaf of groups $\Aut_F(a)$ on
$\cC/X$, where
$$\Aut_F(a)(Y \llra{f} X):= \Aut_{F(Y)}(f^*a).$$

Note also that if $F \ra G$ is an objectwise weak equivalence, then the
induced map of presheaves of homotopy groups is an isomorphism.

\begin{note} If $\cC$ has enough points then a map induces an 
isomorphism on sheaves of homotopy groups if and only if it induces an
isomorphism on the stalks of the sheaves of homotopy groups, which 
is equivalent to inducing an weak equivalence on the stalks.
\end{note}
\begin{reference}[Joyal's Model Structure \cite{Ja}]
There is a left proper, cofibrantly generated,
simplicial model structure on $P(\cC,s\Set)$ where
\begin{itemize}
	\item cofibrations are the maps which are objectwise cofibrations,
	\item weak equivalences are the maps which are isomorphisms on
		sheaves of homotopy groups,
	\item fibrations are the maps with the right lifting property with
		respect to the trivial cofibrations.
\end{itemize}
\end{reference}
The Joyal model category will be denoted by $P(\cC,s\Set)_J$.
Note that in the $S^2$-nullification of $P(\cC,s\Set)_J$ the 
weak equivalences are the maps which induce isomorphisms on
the $sh \pi_0$ and  $sh \pi_1$.

\begin{thm}
\label{AuxMCG}
There is a Quillen equivalence between the model categories
$P(\cC, \Gr)_L$ and the $S^2$-nullification of $P(\cC,s\Set)_J$
given by the adjoint pair $(\fpi, N)$ and the identity 
adjunction of $P(\cC,s\Set)$. 
\end{thm}
\begin{proof}
The result follows from localizing the Quillen equivalence in 
Theorem 1.2 of \cite{DHI} combined with the application of
Corollary A.9 of \cite{DHI}.  
\end{proof}

\begin{cor}
The weak equivalences in $P(\cC, \Gr)_L$ are 
the image under $\fpi$ of those in $(S^2)^{-1}P(\cC,s\Set)_J$.
In particular, a morphism $f \in  P(\cC, \Gr)$ is a local weak 
equivalence if and only if it induces an isomorphism on sheaves of 
homotopy groups.
\end{cor}

\subsection{Characterization of Local Weak Equivalences}

The following definition is the restriction to groupoids of the 
local lifting conditions of 
Section $3$ of \cite{DHI} for simplicial sets.

\begin{defn} \label{local-equivalence-conditions}
A map $F \llra{\phi} G \in P(\cC,\Gr)$ is said to satisfy
the \emph{local lifting conditions} if: 
\begin{enumerate}
\item
Given a commutative square
\[
\parbox[c]{\fill}{
\xymatrix{ \emptyset \ar[r] \ar[d] & F(X) \ar[d] \\ \star \ar[r] & G(X)}}
\hspace{.4cm}
\parbox[c]{\fill}{ \mbox{ $\Rightarrow$ $\exists$ cover $U \ra X$,} }
\hspace{3.3cm}
\parbox[c]{\fill}{
\xymatrix{\star \ar@/^3ex/@{-->}[rrr] \ar[d] & \emptyset \ar[l] \ar[r]
\ar[d] &  F(X) \ar[d] \ar[r] & F(U) \ar[d] \\
\Delta^1 \ar@/_3ex/@{-->}[rrr] & \star \ar[l] \ar[r] & G(X)\ar[r] &
G(U).} }
\]
\\
\item
For $A \ra B$, one of the generating cofibrations 
$\partial \Delta^1 \ra \Delta^1, \text{  } B\Z \ra \star,$ 
given a commutative square
\[
\parbox[c]{\fill}{
\xymatrix{ A \ar[r] \ar[d] & F(X) \ar[d] \\ B \ar[r] & G(X)}}
\hspace{.4cm}
\parbox[c]{\fill}{ \mbox{ $\Rightarrow$ $\exists$ cover $U \ra X$,} }
\hspace{3.3cm}
\parbox[c]{\fill}{
\xymatrix{A  \ar[r] \ar[d] &  F(X) \ar[d] \ar[r] & F(U) \ar[d] \\
B \ar[r] \ar@{-->}[rru] & G(X) \ar[r] & G(U).} }
\]
\end{enumerate}
\end{defn}

\begin{thm} \label{llc-thm}
A map $F \llra{\phi} G \in P(\cC,\Gr)$ is an equivalence on sheaves of 
homotopy groups if and only of it satisfies the local lifting conditions.
\end{thm}

\begin{proof}
Recall that for $F$ a presheaf, its sheafification $shF$,
can be constructed by setting
$$shF(X) = \colim (\eq F(U) \Rightarrow F(V))$$
where the colimit is taken over all covers $U \ra X$ and 
$V \ra U \times_X U$.
It follows that if $a\in shF(X)$ then there exists a cover
$U \ra X$ such that $a$ lifts to an element of $F(U)$.
Similarly if $a,b \in F(X)$ have the same images in $shF(X)$
there exists a cover $U \ra X$ so that they have the same image in $F(U)$.
Conversely these two properties are enough to characterize the
sheafification. It follows that the lifting conditions for
$\emptyset \ra \star$ and $\partial \Delta^1 \ra \Delta^1$ are equivalent
to $sh \pi_0 \phi$ being an isomorphism, and the lifting conditions for
$B\Z \ra \star$ and $\star \ra B\Z$ (which is implied by that for 
$\partial \Delta^1 \ra \Delta^1$) are equivalent to 
$sh \Aut_\phi (a)$ being an isomorphism for all $a \in F(X), X \in \cC$.
\end{proof}

The following two Corollaries are straightforward excercises using 
the local lifting conditions.

\begin{cor} In $P(\cC,\Gr)_L$:
\begin{enumerate}
\item The pullback of a weak equivalence by a levelwise fibration
is again a weak equivalence.
\item The pullback of a weak equivalence which is a levelwise fibration
is weak equivalence which is a levelwise fibration.
\end{enumerate}
In particular $P(\cC,\Gr)_L$ is right proper.
\end{cor}

\begin{cor} \label{lwe-s}
If $F \llra{\phi} G \in P(\cC,\Gr)$ is a local weak equivalence and $F$ is 
a stack then $\phi$ is a levelwise weak equivalence.
\end{cor}

\subsection{Local Model Category Structure on Sheaves of Groupoids}

\begin{prop} \label{comparison}
There exists a local model category structure on $Sh(\cC,\Gr)$, denoted
$Sh(\cC,\Gr)_L$, in which a morphism $f$ is a weak equivlance (resp. fibration)
if and only if $i(f)$ is a weak equivalence (resp. fibration) in $P(\cC,\Gr)_L$.  
Furthurmore, the adjoint pair
$$\xymatrix{P(\cC,\Gr)_L \ar@/_/[rr]_{sh} & & Sh(\cC,\Gr)_L \ar@/_/[ll]_{i} }$$
induce Quillen equivalences between the local model structures.
\end{prop}

\begin{proof}
To see that the model structure $Sh(\cC,\Gr)_L$ is well defined it suffices 
to show that given a generating trivial cofibration $f$ in 
$P(\cC,\Gr)_L$, that $sh(f)$ is a weak equivalence, and that the pushout of 
$sh(f)$ along any morphism in $Sh(\cC,\Gr)$ is still a weak equivalence. 
Both of these statments follow since the natural transformation 
$F \ra i (sh F)$ satisfies the local
lifting conditions, and so is a weak equivalence in $P(\cC,\Gr)_L$.
This also implies that $(sh,i)$ is a Quillen equivalence.
\end{proof}

\begin{cor} \label{we-sheaves}
A morphism $X\llra{f} Y \in Sh(\cC,\Gr)_L$ is a weak equivalence if and
only if it is objectwise full and faithful, and satisfies 
\ref{local-equivalence-conditions}$(1)$.
\end{cor}

The following is a consequence of Corollary \ref{lwe-s}:
\begin{cor}
Let $F \in P(\cC,\Gr)$ be a stack, then $sh(F)$ is also a stack.
\end{cor}

\appendix 

\section{Limits and Colimits in $\Gr/\cC$}

The goal of this section is to prove the following theorem:

\begin{thm}\label{complete}
Categories fibered in groupoids over $\cC$ are closed
under small limits and colimits.
\end{thm}
In order to prove this, we will need a few preliminaries.
\begin{defn}\label{defPFG} $F:\cE \to \cC \in \Cat/\cC$ is pre-fibered 
in groupoids if
\begin{enumerate}
\item Given $f:Y \to X \in \cC$ and $X' \in \cE$ such that $F(X')=X$,
there exists $f' \in \cE$, with target $X'$, such that $F(f')=f$.
\item Given a diagram in $\cE$, over the commutative diagram in $\cC$,
	$$\xymatrix{& Y' \ar[d]^{f'} & \ar@2{->}[r]^F & &
		&  Y \ar[dl]_h \ar[d]^f \\
		Z' \ar[r]^{g'}  & X'  & \ar@2{->}[r]^F & & Z \ar[r]^g & X, }$$
      with $F(f')=f, F(g')=g$, there exists $h'$ such that
      $g' \circ h' = f'$ and $F(h')=h$. Moreover, given two such
      maps $h_1',h_2'$, there exists an automorphism $\phi \in Aut_\cE(Y')$
      such that $F(\phi)=id_Y$ and $h_1' \circ \phi = h_2'$.
\end{enumerate}
\end{defn}

Thus, the difference between fibered and pre-fibered
is that categories which are pre-fibered in groupoids
satisfy a weaker condition than the uniqueness in Condition $(2)$ 
of Definition \ref{defFG}.

\begin{prop}
\label{auxilliary}
Let $I$ be a small category, and $F:I \lra \Gr/\cC$, a diagram.
Then the colimit of the diagram $F$ in $\cCat/\cC$ is pre-fibered in groupoids.
\end{prop}

\begin{proof}
The coproduct in $\cCat/\cC$ of a set of objects in $\Gr/\cC$
is again in $\Gr/\cC$ so it suffices to consider the case of
a coequalizer diagram.
Consider the diagram
$$\xymatrix{ \cR \ar@2{->}[r]^{F_1}_{F_2} \ar[dr] & \cE \ar[r] \ar[d] &
\bar{\cE} \ar[dl] \\ & \cC &}$$
where $F_1,F_2$ are maps in $\Gr/\cC$ and $\bar{\cE}$ is the coequalizer of
the two arrows in $\cCat$.
Recall that the coequalizer in $\cCat$ has objects the coequalizer of
the sets of objects, and morphisms the formal compositions of the
coequalizer of the morphisms, modulo the relations given by composition
in $\cE$. Thus the map $\bar{\cE} \lra \cC$ clearly satisfies Condition $(1)$
of Definition \ref{defPFG}.  The proof that Condition $(2)$ holds follows bellow:

\begin{lemma} \label{factorization} 
Let $X'$ and $Y'$ be in $\cE$, and suppose 
$\bar{f}: [X'] \ra [Y']$ is a map between their images in $\bar{\cE}$.
Write $X=F(X')$ and $Y=F(Y')$ in $\cC$.  Then there is a sequence of
objects $R_1, R_2, \dots R_n \in \cR$ with $F(R_i)= X$ and maps
in $\cE$: 
$$\xymatrix{ X' \ar[r] &  F_1(R_1) & F_2(R_1) \ar[dl] & \\ 
& F_1(R_2) & F_2(R_2) \ar[dl] & \\ & F_1(R_3) & F_2(R_3) \ar[dl] & \\ 
& \dots & F_2(R_{n-1}) \ar[dl] & \\ & F_1(R_n) & F_2(R_n) \ar[r] & Y'}$$
with all but the last covering the identity map of $X$ in $\cC$, such that 
$\bar{f}$ is the composite of images \\
$[X'] \ra [F_1(R_1)]=[F_2(R_1)] \ra [F_1(R_2)]=[F_2(R_2)] \ra [F_1(R_3)] \ra $
\hspace{.1cm} $\dots$  

\hspace{4.5cm}  $[F_2(R_{n-1})] \ra  [F_1(R_n)]=[ F_2(R_n)] \ra [Y']$
in $\bar{\cE}$.
\end{lemma}
\begin{proof}  Use properties of groupoids over $\cC$ to simplify the
expression of a map in $\bar{\cE}$ as formal composition of maps in $\cE$,
modulo the relations given by $\cR$.
\end{proof}
Note that the maps $X' \ra F_1(R_1),$ \hspace{.05cm} 
$F_2(R_1) \ra F_1(R_2),$ \hspace{.05cm} $F_2(R_2) \ra F_1(R_3),$ \hspace{.05cm}
$\dots$ \hspace{.05cm} $F_2(R_{n-1}) \ra  F_1(R_n)$, 
are isomorphisms in $\cE$ because they cover the identity of $X$ in $\cC$.
Also note that $[F_2(R_n)] \ra [Y]$ covers the same maps as $\bar{f}$.

\begin{cor} 
Suppose $\bar{f}$ and $\bar{g}$ are maps $[X] \ra [Y]$
in $\bar{\cE}$ covering the same map $f$ in $\cC$.  
Then there is an automorphism
$\phi$ of $[X']$ such that $ \bar{f} \circ \phi = \bar{g}$.
\end{cor}

\begin{proof}
Factor $\bar{f}$ and $\bar{g}$ as in Lemma \ref{factorization}:\\
$[X'] \ra [F_1(R_1)]=[F_2(R_1)] \ra [F_1(R_2)]=[F_2(R_2)] \ra [F_1(R_3)] \ra $
\hspace{.1cm} $\dots$  

      \hspace{4.5cm} $\dots$ \hspace{.05cm} 
$[F_2(R_{n-1})] \ra  [F_1(R_n)]=[ F_2(R_n)] \ra [Y'],$ \\
$[X'] \ra [F_1(S_1)]=[ F_2(S_1)] \ra [F_1(S_2)]=[ F_2(S_2)] \ra [F_1(S_3)] \ra$

      \hspace{5cm} $\dots$ \hspace{.05cm} 
$[F_2(S_{k-1})] \ra  [F_1(S_k)]=[ F_2(S_k)] \ra [Y'].$\\
Consider the diagram in $\cE$, over the diagram in $\cC$, 
$$\xymatrix{ & F_2(S_k) \ar[d]  & \ar@2{->}[r]^F & & 
& X \ar[d]^f \\ 
F_2(R_n) \ar[r] & Y' &  \ar@2{->}[r]^F & & X \ar[r]^f & Y.}$$
As $\cE$ is fibered in groupoids over $\cC$ there is a unique isomorphism
$F_2(S_k) \llra{h} F_2(R_n)$ covering the identity of $X$ so that 
$$\xymatrix{ & F_2(S_k) \ar[dl]_h  \ar[d] \\ F_2(R_n) \ar[r] & Y'}$$
commutes in $\cE$.
Let $\phi$ be the composite \\
$[X'] \ra [F_1(S_1)]=[F_2(S_1)] \ra$ \hspace{.05cm} 
$\dots$ \hspace{.05cm} $[F_1(S_k)]=[F_2(S_k)] \llra{\bar{h}} 
[F_2(R_n)]=  $ 

\hspace{.5cm} $=[F_1(R_n)] \ra [F_2(R_{n-1})] \ra$ \hspace{.05cm} $\dots$
\hspace{.05cm} $[F_1(R_2)] \ra [F_2(R_1)]=[F_1(R_1)] \ra [X']$ \\
where the second set of maps are the inverses of the isomorphisms
in the factorization of $\bar{f}$.
\end{proof}

\begin{prop} The coequalizer $\bar{\cE}$ satisfies Condition $(2)$ of 
Definition \ref{defPFG}.
\end{prop}
\begin{proof}
Given a diagram in $\bar{\cE}$, over the commutative diagram in $\cC$,
$$\xymatrix{& \bar{Y'} \ar[d]^{\bar{f'}} & \ar@2{->}[r]^F & & &  
Y \ar[dl]_h \ar[d]^f \\ \bar{Z'} \ar[r]^{\bar{g'}} & \bar{X'}  
& \ar@2{->}[r]^F & & Z \ar[r]^g & X, }$$
with $F(\bar{f'})=f, F(\bar{g'})=g$, factor $\bar{f'}$ and $\bar{g'}$ as
in Lemma \ref{factorization} \\
$[Y'] \ra [F_1(R_1)]=[ F_2(R_1)] \ra [F_1(R_2)]=[ F_2(R_2)] \ra [F_1(R_3)]\ra$

\hspace{4cm} $\dots$ \hspace{.05cm} 
$[F_2(R_{n-1})]\ra [F_1(R_n)]=[F_2(R_n)]\ra [X'],$ \\
$[Z'] \ra [F_1(S_1)]=[ F_2(S_1)] \ra [F_1(S_2)]=[ F_2(S_2)] \ra [F_1(S_3)]$ 

\hspace{4.2cm} $\dots$ \hspace{.05cm} 
$[F_2(S_{k-1})] \ra  [F_1(S_k)]=[ F_2(S_k)] \ra [X'].$ \\
Then there is a unique lift of $X \llra{f} Z$ in the diagram in $\cE$,
over the commutative diagram in $\cC$,
$$\xymatrix{ & F_2(R_n) \ar@{-->}[dl]_u \ar[d] &  \ar@2{->}[r]^F & & &  
Y \ar[dl]_h \ar[d]^f  \\  F_2(S_k) \ar[r] & X' & \ar@2{->}[r]^F & &  Z \ar[r]^f & X.}$$
The map 
$[Y'] \ra [F_1(R_1)]=[F_2(R_1)] \ra$ \hspace{.05cm} $\dots$ \hspace{.05cm}  
$[F_1(R_n)]=[F_2(R_n)] \llra{\bar{u}}[F_2(S_k)]$ 

\hspace{2.5cm} $=[F_1(S_k)] \ra [F_2(S_{k-1})]\ra$ \hspace{.05cm} $\dots$ 
\hspace{.05cm} $[F_2(S_1)]=[F_1(S_1)] \ra [Z']$ \\
provides the desired lift of $h$ in $\bar{\cE}$.
\end{proof}

This completes the proof of Proposition \ref{auxilliary}

\end{proof}

\begin{prop} \label{auxilliary2}
Let $\cE \ra \cC$ be pre-fibered in groupoids.
Let $\sim$ be the equivalence relation on $\cE$ generated by
setting $\alpha \sim id$ for the automorphisms $\alpha \in \cE$ which 
satisfy:
\begin{enumerate}
\item  $\alpha$ maps to an identity morphism in $\cC$,
\item there exists  $f \in  \cE $ such
that $f \circ \alpha = f$.
\end{enumerate}
Then $(\cE/\sim) \ra \cC$ is also pre-fibered in groupoids.
\end{prop}

\begin{proof}
The proof follows from the fact that the map 
$\cE \lra (\cE/\sim)$ is surjective on morphisms and bijective on objects.
\end{proof}

\begin{proof}[Proof of Theorem \ref{complete}]
Let $I$ be a small category and $F: I \lra \Gr/\cC$ be a diagram. We
denote by $F'$ the composite $I \llra{F} \Gr/\cC \lra \cCat/\cC$.\\
{\bf Colimits:}
Let $\cE_0$ denote the colimit of $F'$ in $\cCat$.
We will show that the desired colimit of $F$ in $\Gr/\cC$ is the directed 
colimit $\cE$ in $\cCat/\cC$ of the categories $\cE_i$,
where $\cE_i = \cE_{i-1}/ \sim $.
\begin{equation}
\label{coeqcat}
\cE_0  \ra (\cE_0/\sim) = \cE_1 \lra (\cE_1/\sim)=\cE_2 \lra
\cdots  \end{equation}
Propositions \ref{auxilliary} and \ref{auxilliary2} imply that
Condition $(1)$ and the existence part in Condition $(2)$ of Definition
\ref{defFG} are still satisfied by $\cE$.

To show the uniqueness part in Condition $(2)$, suppose given a
commutative diagram in $\cE$:
$$\xymatrix{ Y \ar[r]^f \ar@2{->}[d]^{h_1}_{h_2} & X \\ Z  \ar[ur]^g}$$
such that $h_1$ and $h_2$ project to the same map in $\cC$. 
We can pick lifts $h_1'$ and $h_2'$ of $h_1$ and $h_2$ in some $\cE_i$.
They also project to the same map in $\cC$ and as $\cE_i$ is pre-fibered in
groupoids by Proposition \ref{auxilliary}, there is an automorphism 
$\alpha$ of $Y$ in $\cE_i$ mapping to an identity in $\cC$ such that 
$h_2' \circ \alpha = h_1'$.
It follows that  $h_1' = h_2' \in \cE_{i+1}$
and so $h_1$ and $h_2$ agree in $\cE$.

To show that $\cE$ is the 
colimit in $\Gr/\cC$, observe that if $\cF$ is pre-fibered in groupoids, 
and $\cE'$ is fibered in groupoids, then any map $\cF \ra \cE'\in \cCat/\cC$, 
factors uniquely through $\cF/\sim$.

{\bf Limits:}
Consider the inverse limit of our diagram in $\cCat/\cC$.  
The objects and morphisms of $\lim F'$ are the
inverse limits of the sets of objects and morphisms, so for each object
$X' \in \lim F'$, the category $(\lim F')/X'$, is the inverse limit of
categories $F(i)/X'_i, i \in I$.
It is easy to see that the map $(\lim F' ) /X' \ra \cC/ X$
\begin{itemize}
        \item is a bijection on $\Hom$-sets, since this is the case for each
                of the functors $F(i)/ X'_i \ra \cC/X$,
        \item is not necessarily a surjection on objects, even though
               each of the functors $F(i)/ X_i' \ra \cC/X$ is.
\end{itemize}
Consider the full subcategory of $\lim F'$ with objects
all those $X'$ such that $(\lim F')/X' \ra \cC/X$ is surjective on
objects.  This subcategory is clearly fibered in groupoids and 
satisfies the universal property of the limit.
\end{proof}

\section{Pushouts in $\cCat$}

The goal of this section is to prove:

\begin{prop} \label{B-1}
Let $A, B, C$ be small categories, and $A \llra{i} B$
be a functor which is a monomorphism on objects, and $A \llra{j} C$ a
surjective equivalence of categories.  Then the induced functor
to the pushout in $\cCat$,  $B \ra P:= C \coprod_A B$ is also a
surjective equivalence of categories.
\end{prop}

\begin{proof}
First note that the universal map $B \llra{p} P$ is surjective on objects.
If  $b,b' \in \ob B$, then $p(b)=p(b')$ if and only there
exist $a,a' \in \ob A$ with $i(a)=b, i(a')=b'$ and $j(a)=j(a')$. 
So there is a unique map $a\lra a' \in A$ which maps to the identity
of $j(a)$ and we will call the image of this map in $B$
the canonical map $b \lra b'$. For $b$ not in the image of $A$ the canonical
map $b \lra b$ is defined to be the identity.

It is clear that $p$ induces an isomorphism on components so it
remains to show that $p$ induces an isomorphism
\[ \Hom_B(b,b') \lra \Hom_P(p(b),p(b')). \]

For $\beta,\beta'$ objects of $P$, let
$W(\beta,\beta')$ denote the set of words formed by formal compositions
of morphisms in $B$ and $C$ such that the first map in the word has domain
representing $\beta$, the last map has range representing $\beta'$ and 
consecutive maps have domains and ranges whose images in $P$ agree. Recall that
$\Hom_P(\beta,\beta')$ is the quotient of $W(\beta,\beta')$ by the
equivalence relation generated by the composition in $B$, composition in $C$
and $i(f) \sim  j(f)$ for $f$ a morphism in $A$.

Let $b,b'$ be objects of $B$ and write $\beta=p(b), \beta'=p(b')$. We will
define  functions
$\phi_{b,b'}: W(\beta,\beta') \lra \Hom_B(b,b')$
which are constant on the equivalence classes of
$W(\beta,\beta')$ and so determine functions
$\Hom_P(\beta,\beta') \lra \Hom_B(b,b')$. It will be immediate
from the construction that these are inverse to $p$
and this will complete the proof.

The functions $\phi_{b,b'}$ are defined by induction on the length
of words as follows. Let $w$ be a word of length 1. If $w$ is a
morphism $c \llra{f} c' \in C$ then let $a,a'$ be the unique
objects in $A$ such that $i(a)=b, i(a')=b', j(a)=c, j(a')=c'$ and let
$a \llra{g} a'$ denote
the unique morphism in $A$ such that $j(g)=f$. Define $\phi_{b,b'}(w)=i(g)$.
If $w$ is a morphism $b_1 \llra{f} b_2\in B$ define
$\phi_{b,b'}(w)$ to be the composite $b \lra b_1 \llra{f} b_2 \lra b'$
where the unlabeled arrows are canonical morphisms.

Now suppose $\phi_{b,b'}$ has been defined on words of length $\leq n$
and let $w=w'f$ where $w'$ is a word of length $n$ and $f$ is a morphism
in $B$ or in $C$. Let $b''$ be an arbitrary object of $B$ mapping
to the range of $w'$ and define $\phi_{b,b'}(w)$ as the composite
$b\llra{\phi_{b,b''}(w')} b'' \llra{\phi_{b'',b'}(f)} b'$.
It follows from the construction that the value of $\phi_{b,b'}$
is independent of the choice of $b''$ and that $\phi_{b,b'}$ is
constant on the equivalence classes of $W(\beta,\beta')$.
\end{proof}

\end{document}